\documentclass[11pt,a4paper]{article}
\usepackage[latin1]{inputenc} 
\usepackage[english]{babel}
\usepackage{bm}
\usepackage[T1]{fontenc} 
\usepackage{amsfonts,amssymb,latexsym,epsfig,amsmath}
\usepackage{graphicx}
\usepackage{amsthm}
\usepackage{dsfont}
\usepackage{color}
\usepackage{hyperref}

\theoremstyle{remark}

\numberwithin{equation}{section}
\newcommand{\N}{\hbox{$ I\kern -0.23em N$}}

\newcommand{\Beq}{\begin{equation}}
\newcommand{\Eeq}{\end{equation}}
\newcommand{\BS}{\begin{subequations}}
\newcommand{\ES}{\end{subequations}}
\newcommand{\Beqn}{\begin{equation*}}
\newcommand{\Eeqn}{\end{equation*}}
\newcommand{\Beqa}{\begin{eqnarray}}
\newcommand{\Eeqa}{\end{eqnarray}}
\newcommand{\Beqan}{\begin{eqnarray*}}
\newcommand{\Eeqan}{\end{eqnarray*}}

\newcommand{\tc}{c}

\begin{document}
\title
{Mean flow velocities and mass transport for Equatorially-trapped water waves with an underlying current}

\author{David Henry and Silvia Sastre-G\'omez}

\date{}

\maketitle

\begin{abstract} In this paper we present an analysis of the mean flow velocities, and related mass transport, which are induced by certain Equatorially-trapped water waves. In particular, we examine a recently-derived exact and explicit solution to the geophysical governing equations in the $\beta-$plane approximation at the Equator which incorporates a constant underlying current.  
\end{abstract}
\section{Introduction}
The question of determining the fluid drift induced by the propagation of surface water waves is a fascinating issue and, despite pioneering work on this subject being instigated by Stokes as far back as the mid-$1800$'s, it is still a highly curious and perplexing matter at  even  the most fundamental level. For instance, in the setting of periodic surface gravity water waves  Stokes demonstrated by way of approximations \cite{Sto}  that fluid particles experience a (mean) forward drift to the order of $\epsilon^2$, where $\epsilon$ relates to the wave steepness. This drift is in a mean sense, whereby an average is taken over the wave period, and it is an inherently  nonlinear phenomenon with regard to the order of wave amplitude. The subtleties of these drift properties may be illustrated by considering the classical assumption that for periodic irrotational wave motion it was assumed, at the linear level, that fluid particles follow closed trajectories \cite{John}, whereas according to the Stokes drift phenomenon at the order of expansion $\epsilon^2$ it is implied that at least some particle paths are non-closed. It is noteworthy that, with regard to particle trajectories for periodic irrotational water waves, it was recently proven by various methods that all particle paths throughout the fluid domain are indeed non-closed  for flows induced by a wide-range of gravity (and capillary-gravity) waves, both in the approximate linear regime and for exact solutions of the fully nonlinear governing equations \cite{C2006,Constantin_book,C2012,ConEhVil,ConVil,Hen2006,Hen2008,I-K2008,Lyons}. In recent decades, following the work of Longuet-Higgins, the study of mean drift velocities induced by surface wave motion was placed on a firmer theoretical footing for a broad range of fluid motions  \cite{And,Buh,LS1,LS2}. It was observed that key features of the mean fluid drift velocity, or so-called Stokes' drift velocity, could be characterised in terms of the mean Eulerian flow velocity and the mean Lagrangian flow velocity, whereby: {Lagrange = Euler + Stokes}. In spite of recent progress, determining the mean fluid flow velocities remains a highly complex and intricate issue from both a theoretical, and experimental \cite{Mon,Web}, viewpoint.  

In this paper we present an analysis of the mean flow velocities, and related mass transport, induced by certain Equatorially-trapped water waves. In particular, we examine a recently-derived \cite{Hen2013} exact  solution for the geophysical governing equations in the $\beta-$plane approximation  \cite{ConJohn,Cushman_Roisin,Fedorov_Brown} at the Equator. The form of this solution is explicit in terms of Lagrangian variables, and a benefit inherent in employing the Lagrangian framework is that fluid kinematics may often be described explicitly and with (relative) ease, \cite{Ben,Constantin_edge_2001,Constantin_2001,Con2012,ConJPO,R1,ConstGer,Gerstner,GenHen2014,Henry_2008,Hen2013,HenHsu2015a,HenHsu2015b,I-K,Sil,St}.  A significant complicating factor for the analysis undertaken in this paper, particularly with regard  to determining the mean Eulerian flow velocity and subsequently the  Stokes drift velocity, is the presence of a constant underlying current term in  the solution given in \cite{Hen2013}. This is in spite of the underlying current  assuming a relatively simple manifestation in the Lagrangian formulation of the solution, along the lines of the underlying current terms which were introduced by Mollo-Christensen \cite{Mol} into Gerstner-type solutions in an attempt to model billows and various other complicating effects in both atmospheric and oceanographic situations.  We also note that it is well established that currents play a vital role in Equatorial dynamics \cite{Con2012,R1,ConJohn,Cushman_Roisin,Iz}, and interestingly a transverse Equatorial current can be incorporated into a Gerstner-like exact solution in the Equatorial $f-$plane formulation, cf. \cite{Hen2016}.  The paper is concluded with a brief discussion of some mass-transport properties of these Equatorially trapped waves.

\section{The Equatorially trapped wave solution}
\subsection{Governing equations}
We consider geophysical waves in the Equatorial region, where we assume that the earth is a perfect sphere of radius $R=6378$ km, and work in a  reference
frame rotating with the earth whose origin is fixed at the earth's surface,  with the $\{x,y,z\}$-coordinate frame chosen so that the $x$-axis is  pointing horizontally due east (the zonal direction), the $y$-axis is due  north (meridional direction), and the $z$-axis is pointing vertically upwards and perpendicular to the earth's surface.  The governing equations for geophysical ocean waves are given by 
\begin{subequations}\label{Euler_eq}
\begin{align}
		u_t+uu_x+vu_y+wu_z\,+2\Omega w \cos \Phi
		-2\Omega v\sin\Phi 
		& \displaystyle= -{1\over \rho} P_x\smallskip\\
		v_t\,+uv_x\,+vv_y\,+wv_z\,+2\Omega u \sin\Phi 
		& \displaystyle = -{1\over \rho} P_y\smallskip\\
		w_t+uw_x+vw_y+ww_z-2\Omega u \cos\Phi 
		& \displaystyle = -{1\over \rho} P_z-g,
	\end{align}
\end{subequations}
together with the mass conservation equation
\begin{subequations}\label{Gov}
\begin{equation}\label{mc}
	\rho_t+u\rho_x+v\rho_y+w\rho_z=0	
\end{equation}
and the equation of incompressibility
\begin{equation}\label{in}
	u_x+v_y+w_z=0. 
\end{equation}
Here $\Phi$ represents the latitude, $(u,v,w)$ is the fluid velocity,   
$\Omega=73.10^{-6}$~rad$/$s is the (constant) rotational speed of earth \cite{Cushman_Roisin}, $g=9.8$~m/s$^{-2}$ is the 
gravitational constant, $\rho$ is the water density, and $P$ is the pressure. 
We are interested in Equatorial waves, that is, geophysical ocean waves in a region 
which is within $2^o$ latitude of the Equator. Since the latitude is small, 
we may use the approximations $\sin \Phi\approx \Phi$, and $\cos\Phi\approx 1$, and thus linearising the Coriolis force leads to the $\beta$-plane approximation to equations \eqref{Euler_eq} given by
\begin{equation}\label{beta_plane_eq}
\begin{array}{ll}
		u_t+uu_x+vu_y+wu_z+2\Omega w
		-\beta y v 
		& \displaystyle= -{1\over \rho} P_x\smallskip\\
		v_t+\,uv_x\,+vv_y+\,wv_z+\beta y u
		& \displaystyle = -{1\over \rho} P_y\smallskip\\
		w_t+uw_x+vw_y+ww_z-2\Omega u 
		& \displaystyle = -{1\over \rho} P_z-g,
\end{array}
\end{equation}
where $\beta=2\Omega/R=2.28\cdot 10^{-11}$~m$^{-1}$s$^{-1}$. 
The relevant boundary conditions are the kinematic boundary conditions 
\begin{equation}\label{k}
	w=\eta_t+u\eta_x+v\eta_y \mbox{ on } z=\eta(x,y,t),
\end{equation}
\begin{equation}\label{p}
	P=P_{atm} \mbox{ on } z=\eta(x,y,t),
\end{equation}
where $P_{atm}$ is the (constant) atmospheric pressure, and $\eta(x,y,t)$ 
is the free surface. The boundary condition \eqref
{k} states that all the particles in the surface 
will stay in the surface for all time $t$, and the boundary condition \eqref
{p} decouples the water flow from the motion 
of the air above.  Finally, we assume the water to be infinitely deep, with the flow converging rapidly with depth  to a uniform zonal current, that is, 
\begin{equation}\label{lim}
(u,v,w)\to (-c_0,0,0)\;\mbox{ as }\; z\to-\infty. 
\end{equation}
\end{subequations}
The set of equations \eqref{Gov} comprises the governing equations for the $\beta-$plane approximation of geophysical ocean waves with a constant underlying current.
\subsection{Exact solution}
In this section we briefly describe the exact solution of the $\beta$-plane governing 
equations \eqref{Gov} which was presented in \cite{Hen2013}. This solution prescribes a three-dimensional eastward-propagating steady geophysical wave in the presence of a constant underlying current of magnitude $|c_0|$. The wave-like term is periodic in the zonal direction and it has a constant phasespeed  $c>0$. Furthermore, the wave is Equatorially trapped, exhibiting a strong exponential decay away from the Equator.
Equatorially trapped waves which are symmetric  about the Equator and  propagate eastward are known to exist,  and they are regarded as an important factor in a possible explanation of  the El Ni\~no phenomenon (cf. \cite{ConJohn,Cushman_Roisin, Fedorov_Brown}). 
The  solution of \eqref{Gov} we present is formulated in the Lagrangian framework, whereby the evolution in  time of individual fluid particles  is prescribed \cite{Ben}.
In this Lagrangian formulation  the Eulerian coordinates of fluid particles $(x,y,z)$ are expressed as functions of the Lagrangian labelling variables $(q,r,s)\in \left(\mathbb R,(-\infty,r_0),\mathcal I\right)$, and time $t$,   as follows:
\begin{subequations}\label{exact_sol}
\begin{align}
\label{sol1} x&=q-c_0 t-\frac 1k e^{k[r-f(s)]}\sin{[k(q-\tc t)]},
\\
\label{sol2} y&=s,\\
\label{sol3} z&=r+\frac 1k e^{k[r-f(s)]}\cos{[k(q-\tc t)]},
\end{align}
\end{subequations}
where  $r_0<0$ and $k$ is the wavenumber defined by $k=2\pi/L$, and where $L$ is the (fixed) wavelength. For $c_0>0$ the underlying current is adverse, while for $c_0<0$ the current is { following}, and we see below that the sign of the current determines whether $\mathcal I$ is the real line $\mathbb R$ or a finite interval. The function $f(s)$ determines the decay of the particle oscillations in the latitudinal direction away from the equator and it is given by 
\Beq\label{feq}
f(s)=\frac{\tc \beta}{2\gamma}s^2,
\Eeq
 where 
$
\gamma\!:=2\Omega c_0 +g $$\,\,(>\!0)$ is a ``modified gravity'' term and we make the (physically reasonable) assumption that   $c_0>-\frac{g}{2\Omega}$. 
For notational convenience let us choose 
\[
\xi=k\left(r-f(s)\right),\ \theta=k(q-\tc t).\] Then the Jacobian matrix of the transformation \eqref{exact_sol} is given by 
\Beq \label{jac}
\left({\partial (x,y,z)\over \partial (q,s,r)}\right)
=\left(
\begin{array}{ccc}
1-e^{\xi}\cos \theta & 0 & -e^{\xi}\sin \theta 
\\
f_se^{\xi}\sin \theta & 1 & -f_se^{\xi}\cos \theta  
\\
-e^{\xi}\sin \theta & 0 & 1+e^{\xi}\cos \theta 
 \end{array}
\right),
\Eeq
which has the time-independent determinant $1-e^{2\xi}$. Consequently the flow defined by \eqref{exact_sol} is volume preserving, ensuring that \eqref{in} holds in the Eulerian setting \cite{Ben}. 
Since the solution \eqref{exact_sol} is explicit in the Lagrangian formulation, we may immediately
discern some qualitative properties of the physical fluid motion. Indeed, a significant benefit of working in the Lagrangain framework is that the fluid kinematics can often be  described explicitly and with relative ease. In the case above we calculate the velocity field directly from \eqref{exact_sol} to get
\begin{subequations}
\label{exact_velocity}\begin{align} \label{v1}
u(q,r,s;t)=\frac{Dx}{Dt}&={\tc} e^{\xi}\cos{\theta}-c_0, 
\\ \label{v2} 
v(q,r,s;t)=\frac{Dy}{Dt}&=0, 
\\ \label{v3}
w(q,r,s;t)=\frac{Dz}{Dt}&=\tc e^{\xi}\sin{\theta},
\end{align}
\end{subequations}
where $D/Dt$ is the material (or convective) derivative with respect to Eulerian variables.
For fixed latitudes, that is for every fixed $s$, the system \eqref{exact_sol} describes the flow 
beneath a surface wave propagating eastwards at constant 
speed $c$ determined by the dispersion relation \eqref{aa} below.  Additionally, for fixed latitudes the free surface $z=\eta(x,y,t)$ is obtained by setting $r=r_0(s)$ in \eqref{sol3}, where $r_0(s)<r_0$ is the unique solution to
\Beq \label{sol}
\frac{e^{2k[r(s)-\frac{\tc\beta}{2\gamma}s^2]}}{2k}-r(s)+ \frac{c_0\beta}{2\gamma}s^2-\frac{e^{2kr_0}}{2k}+r_0=0,
\Eeq
The existence of a unique solution $r(s)$ to \eqref{sol} for $|s|>0$ is equivalent to the condition
\Beq \label{cond}
\frac{e^{2k[r_0-\frac{\tc\beta}{2\gamma}s^2]}}{2k}+ \frac{c_0\beta}{2\gamma}s^2<\frac{e^{2kr_0}}{2k},
\Eeq cf. \cite{Hen2013} for details. For $c_0\leq 0$, it is easy to see that condition \eqref{cond} holds for all $s\in \mathbb R$. For $c_0>0$, condition \eqref{cond} will hold for restricted values of $s$ on a finite interval $\mathcal I$ which depends on the magnitude $|c_0|$ of the current. For our present purposes we remark that,  given a current with $c_0>0$,  for a unique solution of \eqref{sol} to exist it is necessary that 
\Beq \label{nec}
c_0 <\tc e^{2kr_0},
\Eeq and accordingly \eqref{exact_sol} represents a dynamically possible 
solution of \eqref{Gov}. 
Since $c_0\neq \tc$ (by \eqref{nec}) it follows that the dispersion relation for the wave takes the form 
\Beq \label{aa}
\tc=\frac{\sqrt{\Omega^2+k\gamma}-\Omega}{k}=\frac{\sqrt{\Omega^2+k(2\Omega c_0 +g)}-\Omega}{k}>0.
\Eeq
We remark that if $c_0=\tc$ then the dispersion relation for the wave would take the form $\tc=\sqrt{g/k}$. Hence, in this situation geophysical Coriolis effects have no bearing on the dispersion relation, which instead matches that of the celebrated Gerstner's wave solution \cite{Constantin_2001,Constantin_book,Henry_2008} for deep-water gravity waves. This observation leads us to infer that precluding the case $c_0=c$, as is consistent with condition \eqref{nec}, is natural in the context  of geophysical water waves (cf. \cite{Hen2013} for details on the dispersion relations). Finally, we note that at fixed latitudes $s=s_*$ the crest and trough levels of the wave surface profile are prescribed in terms of the Lagrangian parameters by
\[
z_{\pm}(s_*)=r_0(s_*)\pm{1\over k}e^{k[r_0(s_*)-f(s_*)]}.
\]

\section{Mean velocities and Stokes drift}
In this section we analyse the effect that the constant underlying current has with respect to both the mean Lagrangian and Eulerian flow velocities induced by the exact solution \eqref{exact_sol}. In \cite{ConstGer} it was shown that in the absence of a current, that is for $c_0=0$, the mean Lagrangian velocity is zero and the mean Eulerian velocity flows westwards. Hence, in the absence of the current the Stokes drift (or mean Stokes flow velocity), which which is the difference in the mean Lagrangian and Eulerian velocities \cite{LS1,LS2}, is eastwards. Here we show that the situation is far more complex in the presence of an underlying current, in particular when determining the mean Eulerian velocity. Throughout the following considerations we fix the latitude by setting $s=s_*$.
\subsection{Mean Lagrangian flow velocity}
The mean Lagrangian flow velocity (also known as the { mass-transport} velocity \cite{LS1}) at a point in the fluid domain is the mean velocity  over a wave period of a marked fluid particle which originates at that point.  For the exact solution \eqref{exact_sol} we may calculate the average of the horizontal velocity $u$ in \eqref{v1} over a wave period $T=L/c$ as follows:
\begin{equation}
	\begin{array}{ll}
	\displaystyle \langle u\rangle_L&\displaystyle={1\over T}\int_0^T u(q-
	ct,s,r)dt\smallskip\\
	&\displaystyle={ce^{\xi}\over T}\int_0^T \!\! \!
	\! \cos\left[k(q-ct)\right] dt-{1\over T}\int_0^T \!\! \!\!c_0\,dt =-c_0,
	\end{array}
\end{equation} where we have  used the fact that the first  integral on the left-hand side above vanishes. It is immediately apparent that the mean Lagrangian flow velocity  is either westwards or eastwards, depending on whether the sign of $c_0$ is positive or negative respectively. When $c_0=0$ the mean Lagrangian velocity is zero, which concurs with the result of \cite{ConstGer}, and in this light the form of the mean Lagrangian flow velocity above is not particularly surprising considering the explicit manner in which $c_0$ appears in the expression for the Lagrangian velocity  \eqref{v1}.  We note that the expression for the mean Lagrangian velocity is independent of both the latitude $s$, and the location   in the fluid domain where the fluid parcel originates.
\subsection{Mean Eulerian flow velocity}
When working in the Eulerian setting matters are greatly complicated by the presence of the underlying current. The mean Eulerian flow velocity at a fixed-point in the fluid domain is the Eulerian fluid velocity at that fixed-point averaged over a wave period. In the case of the velocity field \eqref{exact_velocity} the mean Eulerian flow velocity may be computed by taking  the mean over a wave period of the horizontal velocity \eqref{v1} at any fixed-depth beneath the wave trough. Letting $z=z_-(s_*)$ denote the vertical position of the wave trough level, we fix a depth $z=z_0<z_-(s_*)$. This fixed depth $z=z_0$ may be characterised in terms of Lagrangian variables, using \eqref{sol3}, by the equation
\begin{equation}\label{depth}
	z_0=R+{1\over k}e^{ \xi\left(R\right)}\cos\theta,
\end{equation}
where we denote by $r=R(q-ct; s_*,z_0)$  the functional relationship  induced by relation \eqref{depth} between the otherwise independent variables $r$ and $q$, as follows from the implicit function theorem. We note that a consequence of \eqref{depth} is that $R$ is periodic in the $q-$variable, with period $L$. Differentiating \eqref{depth} with respect to $q$ yields
\begin{equation*}
	0=R_q+R_qe^{\xi\left(R(q)\right)}\cos\theta-e^{\xi\left(R(q)\right)}\sin\theta,
\end{equation*}
that is
\begin{equation}\label{alpha_q}
	R_q={e^{\xi}\sin\theta\over 1+e^{\xi}\cos\theta}.
\end{equation}
We note  from \eqref{alpha_q} that $R$ is maximised or minimised with respect to $q$ whenever $\sin \theta=0$, and therefore  for a fixed-depth $z_0$  the maximal and minimal values achieved by $R$ are given implicitly by the relations
\[
z_0=R\pm{1\over k}e^{ \xi\left(R\right)},
\]
where the positive (negative) sign corresponds to the minimal (maximal) value of $R$, respectively. 
	%
	%
To compute the Eulerian mean velocity $\langle u\rangle_E(s_*,z_0)$ at 
latitude $s_*$ and depth $z_0\le z_-(s_*)$ we examine
\begin{equation*}
\begin{array}{ll}
	  c+\langle u\rangle_E(s_*,z_0)&\displaystyle={1\over T}\int_0^T\left[
	  c+u(x-ct,y,z_0)\right]dt.  
\\
%
%
& \displaystyle={1\over L}\int_0^L\left[
	  c+u(x-ct,y,z_0)\right]dx,  
\end{array}
\end{equation*}
which upon transforming, by way of \eqref{exact_sol}, to the labelling variables $(q,s,r)$,  and invoking functional periodicity with respect to the $q-$variable, we get 
\begin{equation*}
	  \begin{array}{ll}  
	  c+\langle u\rangle_E(s_*,z_0)& \displaystyle={1\over L}\int_{0}^{L}\left[c+u(q-ct,s_*,R(q-ct;s_*,z_0))\right]{\partial x\over\partial q}
	  dq.  
	  \end{array}
\end{equation*}
By differentiating $x$ in \eqref{sol1}  with respect to $q$, using \eqref{v1}, and taking into account \eqref{alpha_q}, we obtain 
\begin{equation*}
\begin{array}{ll}  
	  c+\langle u\rangle_E(s,z_0)\!\!\!\!& \displaystyle={1\over L}\int_{0}^{L
	 }\left[c+ce^{\xi\left(R(q)\right)}\cos\theta-c_0\right]\left[1-e^{\xi\left(R(q)\right)}\cos\theta-e^
	  {\xi\left(R(q)\right)}R_q\sin\theta\right]dq
	  \smallskip\\
	  &  \displaystyle=  {1\over L}\int_{0}^{L
	  }\!\!\!c\left(1+e^{\xi\left(R(q)\right)}\cos\theta\right){1-e^{2\xi\left(R(q)\right)}\over 1+
	  e^{\xi\left(R(q)\right)}\cos\theta}dq-{c_0\over L}\!\int_{0}^{L
	  }\!\!\!{1-e^{2\xi\left(R(q)\right)}\over 1+e^{\xi\left(R(q)\right)}\cos\theta}dq
	  \smallskip\\
	  &  \displaystyle=c-{c\over L}\int_{0}^{L} e^{2\xi\left(R(q)\right)}dq
	  -{c_0\over L}\int_{0}^{L}{1-e^{2\xi\left(R(q)\right)}\over 1+e^{\xi\left(R(q)\right)}\cos\theta}
	  dq.
\end{array}
\end{equation*}
Therefore the mean Eulerian velocity is given by the relation 
\begin{equation}\label{mean_Eulerian_velocity}
	\begin{array}{ll} 
		\displaystyle\langle u\rangle_E(s_*,z_0)&\displaystyle =-{c\over L}\int_
		{0}^{L} e^{2\xi\left(R(q)\right)}dq-
		{c_0\over L}	\int_{0}^{L}{1-e^
		{2\xi\left(R(q)\right)}\over 1+e^{\xi\left(R(q)\right)}\cos\left(k\left[
		q-ct\right]\right)}dq.
	\end{array}
\end{equation}
The presence of a non-zero underlying current $c_0$ adds a significant complicating factor to expression \eqref{mean_Eulerian_velocity}, and in particular the sign (and hence direction) of the mean Eulerian velocity is not easily discernible from the above expression in general. Nevertheless,  depending on the size and direction of the current $c_0$, we may obtain estimates which  determine the direction of the mean Eulerian velocity following from  the inequalities 
\begin{equation}\label{boun_c_0_neg}
	\displaystyle 
	\int_{0}^{L}{1-e^{2\xi}\over 1+e^{\xi}}dq \le \int_{0}^{L}{1-e^{2\xi}\over1+e^{\xi}\cos \theta}dq\le 	\int_{0}^{L}{1-e^{2\xi}\over 1-e^{\xi}}dq.
\end{equation}  
\subsubsection{The case $c_0>0$:}
First of all let us study the case when {$c_0$ is positive}, which represents an  underlying adverse current in the Lagrangian variables. The second integral term on the right-hand side of  inequality \eqref{mean_Eulerian_velocity}  satisfies 
\begin{equation}\label{boun_c_0_pos}
	 -{c_0\over L}
	\int_{0}^{L}{1-e^{2\xi}\over 1-e^{\xi}}dq\le -{c_0\over L}
	\int_{0}^{L}{1-e^{2\xi}\over1+e^{\xi}\cos\theta }dq\le -{c_0\over L}
	\int_{0}^{L}{1-e^{2\xi}\over 1+e^{\xi}}dq.
\end{equation}
Since   $0< c_0<ce^{2kr_0}<c$  from  \eqref{nec}, equation \eqref{boun_c_0_pos} yields 
\begin{equation}\label{boun_Eul_vel_c_0_pos1}
\begin{array}{l}
	\displaystyle	
\langle u\rangle_{E} \displaystyle \le -{c\over L}
	\int_{0}^{L}\!\!\!e^{2\xi}dq-{c_0\over L}
	\int_{0}^{L}{1-e^{2\xi}
	\over 1+e^{\xi}}dq\le -{c_0\over L}
	\int_{0}^{L}{1+e^{3\xi}\over 1+e^{\xi}}dq< 0.
	\end{array}
\end{equation} 
Therefore the mean Eulerian flow velocity is westwards for all admissible values of $c_0$ for which \eqref{nec} holds.
To get an idea of the range of the mean Eulerian flow we note that
\begin{equation}\label{boun_Eul_vel_c_0_pos2}
	\begin{array}{l}
	\displaystyle\langle u\rangle_{E}\ge  -{c\over L}
	\int_{0}^{L} e^{2\xi}dq-{c_0\over L}
	\int_{0}^{L}{1-e^{2\xi}
	\over 1-e^{\xi}}dq
\ge -{c\over L}
	\int_{0}^{L}{1-e^{3\xi}
	\over 1-e^{\xi}}dq.
\end{array}
\end{equation} 
Hence, since $\xi\leq kR<kr_0<0$, we see that for all latitudes $s$ and depths $z_0<z_{-}(s)$ the mean Eulerian flow velocity is in the range
\begin{equation}
	\langle u\rangle_E(s,z_0)\in \left(-c{1-e^{3kr_0}\over 1-e^{kr_0}}, 0\right).
\end{equation}
That the mean Eulerian flow is westward for an adverse current is not surprising, since in the absence of the current the mean Eulerian flow is westward (cf. \cite{ConstGer}) and the presence of the adverse current term in \eqref{mean_Eulerian_velocity} merely serves to exacerbate this effect.

\subsubsection{The case $c_0\leq 0$:}
The case when {$c_0$ is nonpositive}, $c_0\leq 0$, represents an underlying following current. In this case the influence that the current has on the mean Eulerian flow in \eqref{mean_Eulerian_velocity} is  complex and difficult to discern, and it is not generally possible to analytically determine  its effect directly from expression \eqref{mean_Eulerian_velocity}. Nonetheless, we can deduce some broad characteristics of the flow by working as follows. 
The mean Eulerian velocity \eqref{mean_Eulerian_velocity} is { westwards}, that is $\langle u\rangle_E(s_*,z_0)< 0$, if 
\begin{align*}
	-{c_0\over L}\int_{0}^{L}{1-e^{2\xi\left(R(q)\right)}\over 1+e^{\xi\left(R(q)\right)}\cos\theta}dq
&\le	
-{c_0\over L}\int_{0}^{L}{1-e^{2\xi\left(R(q)\right)}\over 1-e^{\xi\left(R(q)\right)}}dq
\\
&\le 
-{c_0}\max_{q\in[0,L]} {{1-e^{2\xi\left(R(q)\right)}}\over {1-e^{\xi\left(R(q)\right)}}}
<	
c \min_{q\in[0,L]} e^{2\xi\left(R(q)\right)} 
\le
{c\over L}\int_{0}^{L} e^{2\xi\left(R(q)\right)}dq.	
\end{align*}
These series of inequalities hold, and accordingly  $\langle u\rangle_E(s_*,z_0)<0$, if
\begin{equation}\label{eq:cond1}
{c_0}> -c \min_{q\in[0,L]} {e^{2k(R(q;z_0)-f(s_*))}{\left(1-e^{k(R(q;z_0)-f(s_*))}\right)} \over {1-e^{2k(R(q;z_0)-f(s_*))}}}.
\end{equation} We note  that in the absence of an underlying current, that is when $c_0=0$,  condition \eqref{eq:cond1} always holds and so the resulting mean Eulerian velocity is always in the westerly direction, an observation which accords with \cite{ConstGer}. 
The mean Eulerian flow \eqref{mean_Eulerian_velocity} is { eastwards}, $\langle u\rangle_E(s_*,z_0)>0$, if 
\begin{align*}
	-{c_0\over L}\int_{0}^{L}{1-e^{2\xi\left(R(q)\right)}\over 1+e^{\xi\left(R(q)\right)}\cos\theta}dq
&\ge	
-{c_0\over L}\int_{0}^{L}{1-e^{2\xi\left(R(q)\right)}\over 1+e^{\xi\left(R(q)\right)}}dq
\\
&\ge 
-{c_0}\min_{q\in[0,L]} {{1-e^{2\xi\left(R(q)\right)}}\over {1+e^{\xi\left(R(q)\right)}}}
>	
c \max_{q\in[0,L]} e^{2\xi\left(R(q)\right)} 
\ge
{c\over L}\int_{0}^{L} e^{2\xi\left(R(q)\right)}dq.	
\end{align*}
These inequalities hold, and hence $\langle u\rangle_E(s_*,z_0)>0$, if
\begin{equation}\label{eq:cond2}
{c_0}< -c \max_{q\in[0,L]} {e^{2k(R(q;z_0)-f(s_*))}{\left(1+e^{k(R(q;z_0)-f(s_*))}\right)} \over {1-e^{2k(R(q;z_0)-f(s_*))}}}.
\end{equation}

 \subsection{Stokes drift}
The Stokes drift (or mean Stokes) velocity $U^S(z_0)$ is defined (cf. \cite{And,ConstGer,LS1,LS2,Mon,Sto}) by the relation
\[
\langle u\rangle_L(z_0)=\langle u\rangle_E(z_0)+U^S(z_0).
\] We derive an expression for the Stokes drift by computing
\begin{align*}
\displaystyle U^S=\langle u\rangle_L-\langle u\rangle_E=
{c\over L}\int_{0}^{L} e^{2\xi\left(R(q)\right)}dq+{c_0\over L}	\int_{0}^{L}{1-e^{2\xi\left(R(q)\right)}\over 1+e^{\xi\left(R(q)\right)}\cos\left(k\left[q-ct\right]\right)}dq-c_0.
\end{align*}
For an adverse current, $c_0\geq0$, it follows from \eqref{nec} that
\begin{align*}
\displaystyle U^S=
{1\over L}\int_{0}^{L} \left(ce^{2\xi\left(R(q)\right)}-c_0\right)dq+{c_0\over L}\int_{0}^{L}{1-e^{2\xi\left(R(q)\right)}\over 1+e^{\xi\left(R(q)\right)}\cos\left(k\left[q-ct\right]\right)}dq>0.
\end{align*}
Therefore for $c_0\geq 0$ the Stokes drift is eastwards throughout the fluid domain. In the case a following current, $c_0<0$, the expression for Stokes drift is altogether more complicated and intractable. Nevertheless we remark that, for $c_0<0$,   if the magnitude of the current is such that \eqref{eq:cond2} holds then the Stokes drift must be westwards.


\section{Mass flux}
We conclude with a brief discussion of mass-transport properties of the flow \eqref{exact_velocity}, where we recall that $\langle u\rangle_L$, being the mean velocity of a marked particle, is sometimes called the 
mass-transport velocity. For a non-zero underlying current, $c_0\neq 0$, we intuitively expect the total mass flux below the free-surface wave past a point  $x=x_0$ fixed in Eulerian coordinates  to be infinite. To see this directly we compute the integral 
\begin{equation}\label{mass_flux_a}
	m(x_0-ct,s)=\int_{-\infty}^{\eta(x_0-ct,s)}u(x_0-ct,s,z)dz.
\end{equation}
In order to transform this expression in terms of the Lagrangian labelling variables we work as follows. 
The implicit function theorem ensures that fixing $x=x_0$ in the expression \eqref{sol1},
\[
x_0=q-c_0t-\frac 1k e^{\xi}\sin \theta,
\] induces a functional relationship between $q$ and the variables $r,t$. Accordingly, we denote $q=\gamma(r,t;s_*)$ and differentiate   with respect to $r$, yielding
\[
	0=\gamma_r-e^{\xi}\sin\theta-\gamma_re^{\xi}\cos\theta,
\]
and so
\[
	\gamma_r={e^{\xi}\sin\theta\over 1-e^{\xi}\cos\theta}.
\]
Using this expression we compute
\[
{dz\over dr}={1+e^{\xi}\cos\theta-\gamma_re^{\xi}\sin\theta}={1-e^{2\xi}
	\over 1-e^{\xi}\cos\theta},
\]
which gives  
\begin{equation}\label{mass_flux}
	\begin{array}{rl}
	m(x_0-ct,s)&\displaystyle=\int_{-\infty}^{r_0}(-c_0+ce^{\xi}\cos\theta){dz
	\over dr}dr\smallskip\\
	&\displaystyle=\int_{-\infty}^{r_0}(-c_0+ce^{\xi}\cos\theta){1-e^{2\xi}
	\over 1-e^{\xi}\cos\theta}dr.
	\end{array}
\end{equation}
Since the terms involving $\xi$ decay exponentially as $r\rightarrow -\infty$, it can be easily seen that the expression \eqref{mass_flux} for the total mass-flux at $x=x_0$ is infinite. This is in stark contrast (as we would expect) to the scenario when there is no underlying current, $c_0=0$, since we deduce from \eqref{sol1} that the function $\gamma$ is now $T-$periodic, and furthermore differentiating \eqref{sol1} with respect to $t$ yields
\begin{equation}\label{gamma_t}
	\gamma_t={-ce^{\xi}\cos\theta\over 1-e^{\xi}\cos\theta}.
\end{equation}
From \eqref{mass_flux} and \eqref{gamma_t} we observe that the mass flux 
is given by 
\begin{equation*}
	m(x_0-ct,s)=\int_{-\infty}^{r_0}-\gamma_t(1-e^{2\xi})dr,
\end{equation*}
and since $\gamma$ is $T-$periodic it follows immediately that the average of the mass flux over a period $T$ is zero (cf. \cite{ConstGer} for full details).
In the case where $c_0$ is non-zero we may still deduce some interesting mass-flow properties near the free-surface. If the magnitude of the current $c_0$ is such that 
\begin{equation}\label{finalcond}
|c_0|\le ce^{k(\tilde r(s_*)-f(s_*))},
\end{equation} where the value $\tilde r(s_*)<r_0(s_*)$ denotes some streamline beneath the surface, then the expression
\begin{equation}\label{xy}	\begin{array}{rl}
	\tilde m(x_0-ct,s)&\displaystyle=\int_{\tilde r}^{r_0}(-c_0+ce^{\xi}\cos\theta){1-e^{2\xi}
	\over 1-e^{\xi}\cos\theta}dr,
	\end{array}
\end{equation}
implies that the mass flux between $\tilde r$ and $r_0$ is positive at the crest and negative at the trough. Therefore, for currents sufficiently small that \eqref{finalcond} holds, and in regions close to the surface between $\tilde r$ and $r_0$, at the crest the mass flux \eqref{xy} is forward and at the trough the mass flux goes backward, matching the properties of the flows observed in \cite{ConstGer,LS2}. In the case in which $c_0<-ce^{k(\tilde r(s_*)-f(s_*))}$  the mass-flux between $\tilde r$ and $r_0$ (given by \eqref{xy}) at the  crest would be forward as usual, however at the trough the mass-flux would also be forward, with this  anomalous behaviour due solely to the constant underlying current.

%

\end{document}